\documentclass[12pt,
a4paper]{amsart}
\usepackage[english]{babel}

\usepackage[latin1]{inputenc}
\usepackage{amsmath}
\usepackage{amscd}
\usepackage{amstext}
\usepackage{amsbsy}
\usepackage{amsopn}
\usepackage{amsthm}
\usepackage{amsxtra}
\usepackage{hyperref}
\usepackage{amsfonts}
\usepackage{amssymb}
\usepackage{mathrsfs}
\usepackage{euscript}
\usepackage{array}
\usepackage{stmaryrd}
\usepackage{verbatim}
\usepackage{natbib}

\theoremstyle{plain}
      \newtheorem{theorem}{Theorem}
      \newtheorem{lemma}[theorem]{Lemma}
      \newtheorem{corollary}[theorem]{Corollary}
      \newtheorem{proposition}[theorem]{Proposition}
      
      \theoremstyle{definition}

      \theoremstyle{remark}
      \newtheorem{remark}[theorem]{Remark}

      \theoremstyle{proof}

\newenvironment{preuves}[1][Proof]{\noindent\textsc{#1 :} }{\ \hfill \rule{0.5em}{0.5em}}

\newcounter{Step}
\def\C{\hspace{.17em}\mbox{l\hspace{-.47em}C} }

\def\I{\mathscr{J}}
\def\H{\mathscr{H}}

\def\tdai{\tilde{\delta}_{\alpha,i}}

\def\H{\mathcal{H}}

\begin{document}
\title[Non-$\Gamma$ under finite non-microstates free Fisher information]{A Note about proving non-$\Gamma$ under a finite non-microstates free Fisher information Assumption }
\author[Yoann Dabrowski]{Yoann Dabrowski*\protect\footnotemark[1]}
\address{UCLA, Department of Mathematics, 520 Portola Plaza, LA, CA 90095}
\email{yoann@math.ucla.edu}
\address{Laboratoire d'Informatique Institut Gaspard Monge, Université
  Paris-Est, 5 bd Descartes, Champs-sur-Marne F-77454 Marne-la-Vallée cedex 2 FRANCE}
\date{}
\maketitle\footnotetext[1]{Research supported in part by NSF grant DMS-0555680.\\ \textit{2000 Mathematics Subject Classification~:} 46L54; \textit{Secondary} : 46L57\\ \textit{Key Words~:}Property $\Gamma$; free Probability; free Fisher information; non-microstates free entropy}
\begin{minipage}[c]{.93\linewidth}
{\scriptsize{\advance\baselineskip -0pt \textsc{\scriptsize Abstract.} \normalfont We prove that if $X_{1},...,X_{n} (n\geq
2)$ are selfadjoints in  a $W^{*}$-probability space with finite
non-microstates free Fisher information, then the von Neumann algebra
$W^{*}(X_{1},...,X_{n})$ they generate doesn't have property $\Gamma$ (especially is not amenable). This is an analog of a well-known
result of Voiculescu for microstates free entropy. We also prove factoriality under finite non-microstates entropy.\par} }
\end{minipage}

\begin{center}\section*{\textsc{Introduction}}
\end{center}

In a fundamental series of papers, Voiculescu introduced analogs of entropy and Fisher information in the context of free probability theory. A first microstates free entropy $\chi(X_{1},...X_{n})$ is defined as a normalized limit of the volume of sets of microstates i.e. matricial approximants (in moments) of the n-tuple of self-adjoints $X_{i}$ in a (tracial) $W^{*}$-probability space $M$. The study of this entropy proved useful for the understanding of the von Neumann algebras $W^{*}(X_{1},...,X_{n})\subset M$ generated by $X_{1},...,X_{n}$, e.g. this leads to the absence of Cartan subalgebras \citep{Vo3} and  primeness \citep{Ge} of free group factors or more generally von Neumann algebras generated by tuples with $\chi(X_{1},...,X_{n})>-\infty$. Voiculescu also extended in \citep{Vo3} the well-known fact that free group factors doesn't have property $\Gamma$ to this class of von Neumann algebras generated by tuples with finite microstates free entropy. Starting from a definition recalled later of free Fisher information \citep{Vo5}, Voiculescu also defined a non-microstates free entropy $\chi^{*}(X_{1},...X_{n})$  with up to now less applications to von Neumann algebras. This entropy is known by the fundamental work \citep{BGC} to be greater than the previous microstates entropy, and believed to be equal (at least modulo Connes' embedding conjecture), so that the question naturally arises of proving the above applications to von Neumann algebras for  $\chi^{*}(X_{1},...,X_{n})>-\infty$. For more details, we also refer the reader to the survey \citep{VoS} for a list of properties as well as applications of free entropies in the theory of von Neumann algebras.

The aim of this note is to prove the easiest result in that direction i.e. under the assumption that the free Fisher Information
$\Phi^{*}(X_{1},...,X_{n})<\infty$ (an assumption stronger than $\chi^{*}(X_{1},...,X_{n})>-\infty$ by a logarithmic Sobolev inequality of \citep{Vo5}), we intend to prove that
$W^{*}(X_{1},...,X_{n})$ doesn't have property
$\Gamma$ (especially is not amenable)
(cf. \citep{Vo3} for the corresponding results in the case of microstates free entropy).
Let us note  that this especially implies that for any $X_{1},...,X_{n}$,
in a $W^{*}$-probability space, and $S_{1}$,...,$S_{n}$ free semicircular
elements free with $X_{1},...,X_{n}$, then $W^{*}(X_{1}+tS_{1},...,X_{n}+tS_{n})$
doesn't have property $\Gamma$ (a result not known, to the best of our knowledge, at least when
$W^{*}(X_{1},...,X_{n})$ is not known to satisfy Connes' embedding conjecture into an ultrapower of the hyperfinite $II_{1}$ factor). We will also prove factoriality under finiteness of non-microstates entropy, especially proving the same kind of degenerate convexity as the one of microstates entropy, i.e. all non extremal states have $\chi^{*}(X_{1},...,X_{n})=-\infty$.

More precisely, % we consider $X_{i}$ self-adjoint non-commutative random
%variables in a tracial $W^{*}$-probability space $M$, and we consider the von Neumann algebra they generate
%$W^{*}(X_{1},...,X_{n})\subset M$.
%In that framework, 
let us recall that $\Phi^{*}(X_{1},...,X_{n})$ is
defined (in \citep{Vo5}) thanks to Hilbert-Schmidt-valued derivations, the so called partial free difference quotients $$\delta_{i}:=\partial_{X_{i}:\C\langle X_{1},...,\hat{X_{i}},...,X_{n}\rangle} :C\langle X_{1},...,X_{n}\rangle\rightarrow
HS(L^{2}(W^{*}(X_{1},...,X_{n})))$$ $$\delta_{i}(X_{j}):=\delta_{ij}1\otimes 1$$ $$
HS(L^{2}(W^{*}(X_{1},...,X_{n})))\simeq L^{2}(W^{*}(X_{1},...,X_{n}),\tau)\otimes L^{2}(W^{*}(X_{1},...,X_{n}),\tau).$$

Thanks to a result of Voiculescu, if $\Phi^{*}(X_{1},...,X_{n})=\sum_{i=1}^{n}||\delta_{i}^{*}1\otimes 1||_2^2<\infty$,
these derivations are closable. And, having in mind of proving first factoriality, if an element, say $Z$, of the center of
$W^{*}(X_{1},...,X_{n})$ were in the domain of $\delta_{i}$, we would write
$0=\delta_{i}([Z,X_{j}])=[\delta_{i}(Z),X_{j}]$ for $j\neq i$ thanks to
Leibniz rule and center property. And thus we would obtain that
$\delta_{i}(Z)$, seen as an Hilbert-Schmidt operator, thus a compact
operator, commutes with a diffuse operator, and thus is zero. A free
Poincaré inequality (due to Voiculescu \citep{Vo6} and recalled later) would imply our
result, that is Z is a scalar times the unit of the von Neumann algebra.

At that point, we have thus to remove the domain assumption assumed valid
on the element $Z$ in the center. In the first section, we prove factoriality under a slightly more general assumption for Fisher information relative to a subalgebra $B$. We will then, in the secondsection, using a variant of Free Poincaré inequality and new boundedness results of (unbounded) dual systems, show our main result according to which
$W^{*}(X_{1},...,X_{n})$ does not have property $\Gamma$. Let us mention
that a previous preprint version of this paper used deeply the notion of $L^{2}$-rigidity
introduced in \citep{P1} to get the same result under a supplementary nonamenability assumption. Here, we thus get nonamenability as a byproduct. Moreover, the third section applies the same tools to prove factoriality under finite non-microstates entropy. We also give a corresponding quantitative inequality in terms of one variant of non-microstates free entropy dimension.

\bigskip

\section{\textsc{Factoriality under finite Fisher Information}}

Let us fix some notations (close to those of \citep{P1}).
We consider $M$ a finite von Neumann algebra with normal faithful tracial
state $\tau$, and $\H$ a $M-M$-bimodule. $D(\delta)$ a weakly dense
*-subalgebra of $M$. 
We suppose here that $\delta:D(\delta)\rightarrow\H$ is a real closable
derivation (real means $\langle \delta(x),y\delta(z) \rangle = \langle \delta(z^{*})y^{*},\delta(x^{*})\rangle$).
$\Delta=\delta^{*}\bar{\delta}$ the corresponding generator of a completely
Dirichlet form, as proved in \citep{S3} (see this paper for the non-commutative
definition of a Dirichlet form, here the Dirichlet form is
$\mathcal{E}(x)=\langle\delta(x),\delta(x)\rangle,
D(\mathcal{E})=D(\Delta^{1/2})$, completely means that $\Delta\otimes I_{n}$
is also the generator of a Dirichlet form on $\mathrm{M}_{n}(M)$).
 Let us introduce a deformation of
resolvent maps (a multiple of a so-called strongly continuous contraction
resolvent, cf e.g. \citep{MaR} for the terminology) $\eta_{\alpha}=\alpha(\alpha + \Delta)^{-1}$, which are
 unital, tracial ($\tau\circ\eta_{\alpha}=\tau$), positive, completely positive maps, and moreover
contractions on $L^{2}(M,\tau)$ and normal contractions on $M$, such that
$||x-\eta_{\alpha}(x)||\leq 2||x||$  and
$||x-\eta_{\alpha}(x)||_{2}\rightarrow_{\alpha\rightarrow \infty}
0$ (as recalled e.g.
in Prop 2.5 of \citep{CiS}). We will also consider $\phi_{t}=e^{-t\Delta}$ the
semigroup of generator $-\Delta$.
 Let us recall two
 relations of  the
resolvent maps (see \citep{MaR} for the first and 
\citep{P1} for the second,
 the integrals are understood as  pointwise
Riemann integral)~:

$$\forall\alpha>0, \eta_{\alpha}=\alpha \int_{0}^{\infty}e^{-\alpha t} \phi_{t} dt.$$

$$\forall\alpha>0, \zeta_{\alpha}:=\eta_{\alpha}^{1/2}=\pi^{-1}\int_{0}^{\infty}\frac{t^{-1/2}}{1+t}\eta_{\alpha(1+t)/t}dt$$

%We then use the notation, following J. Peterson, $\zeta_{\alpha}=\eta_{\alpha}^{1/2}$ and $\tilde{\delta_{\alpha}}=\alpha^{-1/2}\delta\circ\zeta_{\alpha}$.
The point is that $Range(\eta_{\alpha})=D(\Delta)\subset D(\bar{\delta})$
and $Range(\eta_{\alpha}^{1/2})=D(\Delta^{1/2})= D(\bar{\delta})$
and  so that  $\bar{\delta}\circ\zeta_{\alpha}$ is bounded (remark that this way to
precompose with $\eta_{\alpha}^{1/2}$ to extend a map to the whole space is
usual in classical Dirichlet form theory (especially in the relation with
Malliavin calculus), in that way, for instance, the gradient operator of
Malliavin calculus is extended to a distribution valued operator (after
post-composition with another operator)).

We now prove the first theorem of that note :

\begin{theorem}\label{Fact}
Let $(M,\tau)$ a tracial $W^{*}$-probability space (i.e. M a von Neumann
algebra with $\tau$ a faithful tracial normal state). Let $(X_{1},...X_{n})$
a n-tuple (of self-adjoints) such that the microstates free Fisher information
$\Phi^{*}(X_{1},...,X_{n})<\infty$, then $W^{*}(X_{1},...,X_{n})$ is
a factor.
\end{theorem}

\begin{preuves}
Let
$\delta_{i}=\partial_{X_{i}:\C\langle X_{1},...,\hat{X_{i}},...,X_{n}\rangle}$
following the notation of Voiculescu for the non-commutative difference
quotient. We see $\delta_{i}:C\langle X_{1},...,X_{n}\rangle\rightarrow
HS(L^{2}(W^{*}(X_{1},...,X_{n})))\simeq L^{2}(W^{*}(X_{1},...,X_{n}),\tau)\otimes L^{2}(W^{*}(X_{1},...,X_{n}),\tau)$.

First, thanks to a result of Voiculescu, $\Phi^{*}(X_{1},...,X_{n})<\infty$
implies that all the derivations $\delta_{i}$ are closable as unbounded operators $L^{2}(W^{*}(X_{1},...,X_{n}),\tau)\rightarrow
HS(L^{2}(W^{*}(X_{1},...,X_{n})))$ and they are even real closable derivations.

But let us now fix $i$ and consider $Y\in\C\langle X_{1},...,\hat{X_{i}},...,X_{n}\rangle$.
By definition, we have $\delta_{i}Y=0$, so that if
$\Delta_{i}=\delta_{i}^{*}\bar{\delta_{i}}$, we have especially
$\Delta_{i}Y=0$ (and $Y\in D(\Delta_{i})$).

Using a complete positivity argument (or an easy differential equation argument on the corresponding semigroup) one can easily show that $\zeta_{\alpha,i}(ZY)=\zeta_{\alpha,i}(Z)Y$ and
$\zeta_{\alpha,i}(YZ)=Y\zeta_{\alpha,i}(Z)$. Thus, 
$\zeta_{\alpha,i}([Z,Y])=[\zeta_{\alpha,i}(Z),Y]$, and if we note
$\tdai=\alpha^{-1/2}\delta_{i}\circ\zeta_{\alpha}$ (a
bounded map as already noted), we have, using Leibniz rule and
$\delta(Y)=0$: 
$$\tdai([Z,Y])=[\tdai(Z),Y]$$

Consequently, if $Z$ is in the center of $W^{*}(X_{1},...,X_{n})$, 
we have proved $[\tdai(Z),Y]=0$. But now, if $Y=X_{j}$ ($j\neq i$), $Y$
is diffuse (inasmuch as $\Phi^{*}(X_{1},...,X_{n})<\infty$ implies
$\chi^{*}(X_{1})+...+\chi^{*}(X_{n})\geq\chi^{*}(X_{1},...,X_{n})>-\infty$,
and if $\xi \in L^{2}(W^{*}(X_{1},...,X_{n}))$ were an eigenvector of
$X_{j}$ with eigenvalue $\lambda$, the projector on $\xi$ in $B(L^{2}(W^{*}(X_{1},...,X_{n}))$ were
not zero, implying the spectral projection $1_{X_{j}=\lambda}$ to be not
zero, and by faithfulness $\tau(1_{X_{j}=\lambda})\neq 0$ a contradiction,
since $\chi^{*}(X_{j})> -\infty$ implies that the distribution of $X_{j}$
has no point masses).

But now, an Hilbert-Schmidt (thus compact) operator commuting with a
diffuse one is zero (using the spectral theorem for compact operators, the
diffuse one should have an eigenvector !).

We have eventually proved $\tdai(Z)=0$ for all $i$ (and all $\alpha>0$) as
soon as $Z$ is in the center of $W^{*}(X_{1},...,X_{n})$ and
thus, by closability, knowing  $||Z-\zeta_{\alpha,i}(Z)||_{2}\rightarrow_{\alpha\rightarrow \infty}
0$, we obtain the fact that $Z\in D(\bar{\delta_{i}})$ and $\bar{\delta_{i}}(Z)=0$.
Then, we conclude with the following lemma, due to Voiculescu (unpublished \citep{Vo6})
\end{preuves}

\bigskip

\begin{lemma}[Free Poincaré inequality \citep{Vo6}]\label{FPI}
Consider $\delta_{i}$ the partial free difference quotient with respect to
$X_{1},...,X_{n}$, and $Y$ a self-adjoint variable in the domain of all the
operators $\bar{\delta_{i}}$ (as unbounded operators
$L^{2}(W^{*}(X_{1},...X_{n}))\rightarrow HS(L^{2}(W^{*}(X_{1},...X_{n})))$), then, there exists a positive constant $C$ depending on the $X_{i}$ but not on $Y$
such that :$$C\sum_{j=1}^{n} ||\bar{\delta_{j}}Y||_{HS}\geq ||Y-\tau(Y)||_{2}$$
\end{lemma}

We refer the reader to \citep{Vo6} for a proof, but we note that the key tool is the following remark, that 
for a polynomial $Y=P(X_{1},...,X_{n})\in \C\langle X_{1},...,X_{n}\rangle$, we verify immediately by linearity and monomial case that :
\begin{equation}\label{FPE}\sum_{j=1}^{n}((\delta_{j}P)(X_{j}\otimes 1) - (1\otimes
X_{j})(\delta_{j}P))=P\otimes 1 - 1\otimes P.\end{equation}

\bigskip

%With the notations of the proof, note for further use that equation (\ref{FPE}) can also be extended, we have equation
%(\ref{FPE}) for $Z_{k}$, and
%$|\left(\delta_{i}(Z_{k})-\bar{\delta_{i}}Y\right)(X_{i}\otimes
%1)|_{2}\leq ||\delta_{i}(Z_{k})-\bar{\delta_{i}}Y||_{2}||X_{i}||\rightarrow 0$ %(note that the multiplication by $X_{j}\otimes
%1$ on $L^{2}(W^{*}(X_{1},...,X_{n}),\tau)\otimes
%L^{2}(W^{*}(X_{1},...,X_{n}),\tau)$ is well-defined and bounded thanks to
%the boundedness already remarked on polynomials and thanks to density).

%Likewise, $||Z_{k}\otimes 1-Y\otimes 1||_{2}\rightarrow 0$, and we have
%obtained :

%\begin{lemma}\label{FPL}
%Consider $\delta_{i}$ the partial free difference quotient with respect to
%$X_{1},...,X_{n}$, supposed closable as before, and $Y$ a self-adjoint variable in the domain of all the
%operators $\bar{\delta_{i}}$ (as unbounded operators
%$L^{2}(W^{*}(X_{1},...X_{n}))\rightarrow
%HS(L^{2}(W^{*}(X_{1},...X_{n})))$), then :
%$$\sum_{j=1}^{n}((\bar{\delta_{j}}Y)(X_{j}\otimes 1) - (1\otimes
%X_{j})(\bar{\delta_{j}}Y))=Y\otimes 1 - 1\otimes Y.$$
%\end{lemma}

\begin{remark}
As Jesse Peterson pointed out to us, after reading  an earlier version of this note,
once we have shown $\bar{\delta_{1}}(Z)=0$ (using only commutation with
$X_{2},...,X_{n}$, we can conclude by writing down
$0=\bar{\delta_{1}}([Z,X_{1}])=[\bar{\delta_{1}}(Z),X_{1}]+[Z,1\otimes 1]=[Z,1\otimes
1]$ and conclude taking the $||.||_{2}$ norm. (We have used our original
proof inasmuch as a variant of free Poincaré inequality will be essential
in the next part. But somehow, the following result is the only one not provable
under weakened assumptions in what follows. The counterpart of the powerfulness of free Poincaré like technique being its non applicability in the case $\Phi^{*}(X:B)<\infty$, up to now). We have thus also proved the following result :
\end{remark}

\begin{theorem}
Let $X$ a selfadjoint element in $M$ a tracial $W^{*}$-probability space, and $B$ a subalgebra of M,
algebraically free with $X$ and containing a diffuse element. Suppose
moreover that the free Fisher information of X relative to B :
$\Phi^{*}(X:B)<\infty$, then $W^{*}(X,B)$ is a factor.
\end{theorem}

Let us end this part with a corollary. Consider, e.g. as in \citep{Hiai}, the full
(universal) free-product $C^{*}$-algebra
$C([-R,R])^{\star N}$ and
note $\mathcal{T}$ the space of tracial states on this $C^{*}$-algebra. It
is (elementarily known to be) a compact convex set for the weak-$*$ topology. It
is moreover known by the reduction theory for von Neumann algebras that
this is a Choquet simplex.  It
is known (see the beginning of the proof of  Th  3.1.18 of \citep{Sakai}
using mainly Prop 3.1.10 and the discussion before the Theorem) that factorial
states (i.e. states for which the bicommutant in the GNS construction give factors) are exactly
extreme points of this convex set. We will show that $\mathcal{T}$ is a
Poulsen Simplex (see \citep{Lin}), i.e. a metrizable Choquet Simplex in which
the extreme points form a dense set. Since, as a weak-$*$ compact of the
dual of a separable space, $\mathcal{T}$ is metrizable, we have merely to
prove the last statement about density of theset of  extreme points.

Let us prove this in the following :

\begin{corollary}
The set of tracial states $\mathcal{T}$ on $C=C([-R,R])^{\star N}$ is
a Poulsen Simplex.
\end{corollary}

\begin{preuves}
To conclude the proof, consider thus a tracial state $\tau$ on $C$,
consider $X_{i}$, the function $f_{0}(t)=t$ in the i-th copy of $
C([-R,R])$, the usual self-adjoint generators (as a $C^{*}$-algebra) of
$C$. Let $W$ the weak closure of $\pi_{\tau}(C)$, the GNS construction
associated to $\tau$, we always note $\tau$ the associated (faithful normal
tracial) state on $W$. We can consider the von Neumann algebra $M$ generated by
$W$ and a free family of semicircular elements $M=W^{*}(W,S_{i})$, and get
another faithful tracial state on $M$ (the last one noted $\tau$,
cf. \citep{CRM} for faithfulness). Consider
$Y_{i,t}=\frac{R(X_{i}+t S_{i})}{R+2t}$. By corollary 3.9 in \citep{Vo5},
$\Phi^{*}(Y_{1,t},...,Y_{N,t})<\infty$ and thus by our theorem \ref{Fact},
$W^{*}(Y_{1,t},...,Y_{N,t})\subset M$ is a factor. But since $||Y_{i,t}||\leq R$,
we have a *-homomorphism $C\rightarrow M$ sending $X_{i}$ to $Y_{i,t}$ (e.g. Prop
2.1 in \citep{Hiai}). This defines by composition with the state on $M$, a
state $\tau_{t}$ on $C$. Since the state considered on $M$ is faithful, the
kernel of the $*$-homomorphism is nothing but the ideal of elements with
$\tau_{t}(Z^{*}Z)=0$ by which we quotient $C$ in the GNS construction for
$\tau_{t}$ on $C$, we thus get a $*$-isomorphism, from this quotient on its
image which preserves the trace, and thus, $L^{2}(W^{*}(Y_{1,t},...,Y_{N,t}),\tau)$ is
isomorphic to $L^{2}(C,\tau_{t})$ (by the induced map). And a step further we get the $*$-isomorphism
between $W^{*}(Y_{1,t},...,Y_{N,t})$ and $\pi_{\tau_{t}}(C)$. Thus,
$\tau_{t}$ is a factorial, thus an extremal tracial state in $\mathcal{T}$. Now,
to get weak-$*$ convergence of $\tau_{1/n}$ to $\tau$ in $\mathcal{T}$, and
thus the claimed density, we
have just to consider convergence on monomials on which a usual trick
shows the concluding inequality
: \begin{align*}\left|\tau(X_{i_{1}}...X_{i_{p}})-\tau_{t}(X_{i_{1}}...X_{i_{p}})\right|&=\left|\tau(X_{i_{1}}...X_{i_{p}})-\tau(Y_{i_{1},t}...Y_{i_{p},t})\right|\\ &\leq R^{p-1}p\sup_{i}(||\frac{R(X_{i}+t S_{i})}{R+2t}-X_{i}||\\ &\leq R^{p-1}p\frac{4Rt}{R+2t}.\end{align*}
\end{preuves}

\section{\textsc{Non-$\Gamma$}}

In the preceding part, we used the semi-group and resolvent maps associated
to a derivation $\delta_{i}=\partial_{X_{i}:\C\langle
  X_{1},...,\hat{X_{i}},...,X_{n}\rangle}$. The drawback is that, if we
have not exactly a commutator equal to zero, as this is the case when we want to prove
non-$\Gamma$, we cannot move the estimate on the commutator (giving an
estimate on the Hilbert-Schmidt operator), to an estimate
on $Z-\tau(Z)$ using something like free Poincaré inequality, inasmuch as we have not the
same resolvent maps for different derivations $\delta_{i}$. We have thus
searched to move (somehow) the preceding reasoning
in case we consider $\delta:=(\delta_{1},...,\delta_{n})$, the resolvent map
associated to it $\eta_{\alpha}(Z)$, and then
$\delta_{i}\circ\eta_{\alpha}(Z)$. This was our approach in a previous version of this paper where we assumed non-amenability and used then a $L^{2}$-rigidity technique to conclude. However, working a little bit more from the following variant of free Poincaré inequality will be much more efficient, enabling us to prove non-$\Gamma$ without any other assumption, and thus proving non-amenability instead of assuming it.

\subsection{Two preliminaries}
First in order to get our inequality, we will use the following general result about derivations in von Neumann algebras, which can be thought of as a ``Kaplansky density theorem for derivations'' (the proof also confirms this), which is of independent interest and really likely known to specialists but for which we have not found any reference.

\begin{proposition}\label{Wdev}
  Let $\delta$ a symmetric derivation defined on a weakly dense $*$-algebra $D(\delta)$ of the tracial $W^{*}$-probability space  $(M,\tau)$, closable as an operator $D(\delta)\subset L^{2}(M,\tau)\rightarrow\H$, $\H$ an involutive  $M-M$ Hilbert $W^{*}$-bimodule (with isometric involution, as usual we assume $\sigma$-weak continuity of both actions). Then the following properties are equivalent :\begin{enumerate}\item[(i)] $D(\bar{\delta})\cap M$ is a $*$-algebra on which $\bar{\delta}|_{D(\bar{\delta})\cap M}$ is a (symmetric) derivation.
\item[(ii)]for any $Z\in D(\bar{\delta})\cap M$, there exists a sequence $Z_{n}\in D(\delta)$ with $||Z_{n}||\leq ||Z||$, $||Z_{n}-Z||_{2},||\delta(Z_{n})-\bar{\delta}(Z)||_{2}\rightarrow 0$. 
\end{enumerate}
\end{proposition}

\begin{preuves}
The fact that $(ii)$ implies $(i)$ is clear since taking $Z_{n}$, $Y_{n}$ for $Z$ and $Y$, as in $(ii)$, $||Z_{n}Y_{n}-ZY||_{2}\rightarrow 0$, and for any $\xi\in \H$, we get successively $||\xi (Y_{n}^{*}-Y^{*})||_{\H}\rightarrow 0$, by coincidence of $L^{2}$ and $\sigma$-$*$-strong topologies on bounded sets in $M$ and $\sigma$-weak (thus $\sigma$-$*$-strong) continuity of the action, and thus $\langle \delta(Z_{n}) ,\xi Y_{n}^{*} \rangle \rightarrow \langle \delta(Z) ,\xi Y^{*} \rangle$, which gives at the end weak convergence of $\delta(Z_{n}Y_{n})$ to $Z\bar{\delta}(Y)+\bar{\delta}(Z)Y$, and by (weak) closability of the graph of $\bar{\delta}$, we get $ZY\in D(\bar{\delta})\cap M$ with the derivation property.

The proof of the converse follows verbatim the proof of Kaplansky density theorem. We may first assume that $\delta:D(\delta)\rightarrow L^{2}(M)$ is closed as a derivation $M\rightarrow L^{2}(M)$, since it is closable and then obtaining $Z_{n}$ in this enlarged $D(\delta)$ is harmless. We can also assume $||Z||\leq 1$. Consider $X=(1+(1-ZZ^{*})^{1/2})^{-1}Z\in M$, then $X\in D(\delta)$ by closability of any closed derivation on a $C^{*}$-algebra by $C^{1}$-functional calculus. Look at $f(x)=2(1+xx^{*})^{-1}x=2x(1+x^{*}x)^{-1}$ so that $Z=f(X)$. Take $X_{n}\in D(\delta)$ converging to $X$ in $L^{2}$ with $||\bar{\delta}(X_{n}-X)||_{2}\rightarrow0$. Then consider $Z_{n}=f(X_{n})$ so that $||Z_{n}||\leq 1$ and $Z_{n}$ converge to $Z$ in $L^{2}$ (even if this is not a consequence of $*$-strong continuity of $f$ since we don't know whether $X_{n}$ converge to $X$ $*$-strongly since it is not bounded as a sequence in $M$, the proof is however standard, like in the proof of Kaplansky density Theorem, see bellow for an example for the derivative of $f$). Since, by hypothesis, $\delta$ and $\bar{\delta}|_{D(\bar{\delta})\cap M}$ are derivations, closed seen as derivation $M\rightarrow \H$, we get $\delta(Z_{n})=2\delta(X_{n})(1+X_{n}^{*}X_{n})^{-1}- 2 X_{n}(1+X_{n}^{*}X_{n})^{-1}\delta(X_{n}^{*}X_{n})(1+X_{n}^{*}X_{n})^{-1}$ and the analog for $\bar{\delta}(Z)$, using appropriate series expansions. Now the boundedness as sequences in $M$ of $X_{n}(1+X_{n}^{*}X_{n})^{-1}X_{n}^{*}$,$(1+X_{n}^{*}X_{n})^{-1}$ and $X_{n}(1+X_{n}^{*}X_{n})^{-1}$ shows that it suffices to show the convergence of $2\bar{\delta}(X)(1+X_{n}^{*}X_{n})^{-1}- 2 X_{n}(1+X_{n}^{*}X_{n})^{-1}(\bar{\delta}(X^{*})X_{n}+X_{n}^{* }\bar{\delta}(X))(1+X_{n}^{*}X_{n})^{-1}$. Likewise, by coincidence of $L^{2}$ and $\sigma$-$*$-strong topologies on bounded sets in $M$ and $\sigma$-$*$-strong continuity of the action, it suffices to show convergence in $L^{2}(M)$ of $X_{n}(1+X_{n}^{*}X_{n})^{-1}X_{n}^{*}$,$(1+X_{n}^{*}X_{n})^{-1}$ and $X_{n}(1+X_{n}^{*}X_{n})^{-1}$. Let us for instance prove the first one, let us write :
\begin{align*}X_{n}&(1+X_{n}^{*}X_{n})^{-1}X_{n}^{*}-X(1+X^{*}X)^{-1}X^{*}=X_{n}(1+X_{n}^{*}X_{n})^{-1}(X_{n}^{*}-X^{*})\\ &+X_{n}(1+X_{n}^{*}X_{n})^{-1}\left[(X^{*}-X_{n}^{*})X+X_{n}^{*}(X-X_{n})\right](1+X^{*}X)^{-1}X^{*}\\ &+(X_{n}-X)(1+X^{*}X)^{-1}X^{*}.\end{align*}
Since for each term, both sides of $X-X_{n}$ or its adjoint are bounded by functional calculus, this concludes.

\end{preuves}

\begin{remark}
Let us note that in order to apply the previous proposition to the free difference quotient in case of finite Fisher information, we prove $(i)$ using Proposition 3.4 in \citep{DL} to get $D(\bar{\delta})\cap M$ is an algebra (knowing that $D(\bar{\delta})$ is the domain of a Dirichlet form, as already recalled, thanks to \citep{S3}) , and then, for instance use the formula for $\delta^{*}$ given by Corollary 4.3 in \citep{Vo5} to show $\bar{\delta}$ is closable as an operator valued in $L^{1}(M\otimes M)$, and we prove there the derivation property, deducing it for the derivation valued in $L^{2}$ as a consequence.
\end{remark}

Second, we recall for the reader convenience some results about bounded and unbounded dual systems in the sense of Voiculescu and Shlyakhtenko respectively. Even if we will not use their results explicitly, this will enable to express some assumptions and results in terms of these standard objects. Let us recall the following result of \citep{Shly} (deduced from Theorem 1 and its proof), $\delta_{i}$ the i-th partial difference quotient as earlier.

\begin{proposition}[\citep{Shly}]
$\delta_{i}^{*}1\otimes 1$ exists (in $L^{2}(W)$, $W=W^{*}(X_{1},..,X_{n}$) if and only if there exists a closable unbounded operator $Y_{i} : L^{2}(W)\rightarrow L^{2}(W)$ with $\C\langle X_{1},...,X_{n}\rangle\subset D(Y_{i})$, $Y_{i}1=0$, $1\in D(Y_{i}^{*})$ such that $[Y_{i},X_{j}]=\delta_{i}(X_{j})$. Moreover, necessarily such a $Y_{i}=1\otimes \tau\circ \delta_{i}$ (or is an extension of it beyond $\C\langle X_{1},...,X_{n}\rangle$).
\end{proposition}

Then corollary 1 of the same paper noticed that $\tilde{Y}_{i}=\frac{1}{2}(Y_{i}-Y_{i}^{*})$ is an anti-symmetric closable dual system. Moreover the proof of Theorem 1 also shows that $Y_{i}^{*}X=X\delta_{i}^{*}1\otimes1 -Y_{i}X$ (also a consequence of Corollary 4.3 in \citep{Vo5} in the free difference quotient case we are interested in here), so that  $\tilde{Y}_{i}(X)=Y_{i}(X)-\frac{1}{2}X\delta_{i}^{*}1\otimes1$ i.e. $\tilde{Y}_{i}1=-\frac{1}{2}\delta_{i}^{*}1\otimes1$. Moreover, it is easily seen that such a $\tilde{Y}_{i}$ gives in inverting the above relation to get a $Y_{i}$ a $Y_{i}$ similar to the one in the previous proposition, we will thus be later interested in bounded dual systems in the sense of Voiculescu verifying this relation for the specific relation they have with the canonical dual system of the previous proposition (for which we will get latter e.g. nice boundedness properties).

\subsection{A mixed Poincaré-non-$\Gamma$ (in)equality}

Our main tool will be a lemma based on the same argument as free Poincaré
inequality. After proving it, we develop several consequences under (more or less) stronger assumptions for further use.

\begin{lemma}\label{FPC}
Let $(M,\tau)$ a tracial $W^{*}$-probability space. Let $(X_{1},...X_{n})$
a n-tuple (of self-adjoints, $n\geq 2$ in order to have a non-trivial result) such that the microstates free Fisher information
$\Phi^{*}(X_{1},...,X_{n})<\infty$. Let $Z \in W^{*}(X_{1},...,X_{n})\cap
D(\bar{\delta})$  ($\delta$ the free difference quotient), then
we have the following equality :
\begin{align*}2&(n-1)||Z-\tau(Z)||_{2}^{2}\\&=\sum_{i=1}^{n} <[Z,X_{i}],[Z,\Delta(X_{i})]>+2 \Re<(\tau\otimes 1-1\otimes\tau)(\bar{\delta_{i}}(Z)),[Z,X_{i}]>.\end{align*}
\end{lemma}

\begin{preuves}It suffices to show the result for $Z\in \C\langle X_{1},...,X_{n}\rangle$ (using Proposition \ref{Wdev}).
Inasmuch as $\delta_{i}$ is a derivation, we have
$\delta_{i}[Z,X_{i}]=[\delta_{i}(Z),X_{i}]+[Z,1\otimes  1]$ and we have
already noticed that :$||[Z,1\otimes  1]||^{2}_{HS}=2||Z-\tau(Z)||_{2}^{2}$.

Everything will be based on the equality on which is based free Poincaré
inequality. Let us compute $||[Z,1\otimes  1]||^{2}_{HS}=<\delta_{i}[Z,X_{i}]-[\delta_{i}(Z),X_{i}],[Z,1\otimes  1]>$:
\begin{align*}<[\delta_{i}(Z),X_{i}],&[Z,1\otimes  1]>=<\delta_{i}(Z),[[Z,1\otimes  1],X_{i}]>
\\ &=<\delta_{i}(Z),[Z,[1\otimes1,X_{i}]]>-<\delta_{i}(Z),[1\otimes1,[Z,X_{i}]]>
\end{align*}
At that point, we notice that : \begin{align*}[Z,[1\otimes
1,X_{i}]]&=[Z,1\otimes X_{i}- X_{i} \otimes 1]\\
&=Z\otimes X_{i}- ZX_{i} \otimes
1-1\otimes X_{i}Z+X_{i} \otimes Z\\
&=(1\otimes X_{i})[Z,1\otimes
1]- [Z,1\otimes1](X_{i} \otimes 1)
\end{align*}
But now, we have in fact written an ``inner'' commutant of $X_{i}$ and $[Z,1\otimes
1]$ (i.e. a commutant with the action of the von Neumann algebra on
$M\otimes M$ on the side of the tensor product not on the outer side,
remark that the preceding equation is just commutation of the two actions
after writing $[1\otimes
1,X_{i}]$ in terms of an inner commutant).

We will merely now use that the scalar product of Hilbert Schmidt operators is
compatible with this inner commutant (which is nothing but an extension
of  traciality of
$\tau\otimes\tau$ on $M\otimes M$):
\begin{align*}\sum_{i=1}^{n}\langle
\delta_{i}(Z),[Z,[1\otimes 1,X_{i}]]\rangle&=
\langle
\sum_{i=1}^{n}(1\otimes X_{i})\delta_{i}(Z)- \delta_{i}(Z)(X_{i} \otimes 1),[Z,1\otimes
1]\rangle\\ &=\langle
(1\otimes Z- Z\otimes 1),[Z,1\otimes
1]\rangle\\ &=-||[Z,1\otimes
1]||^{2}\end{align*}
We have used the equation (\ref{FPE}) on which is based the proof of free
Poincaré inequality.
Thus, we have obtained :
$$\sum_{i=1}^{n}<[\delta_{i}(Z),X_{i}],[Z,1\otimes  1]>=-||[Z,1\otimes
1]||^{2}-\sum_{i=1}^{n}<\delta_{i}(Z),[1\otimes1,[Z,X_{i}]]>.$$
We have now to compute \begin{align*}<\delta_{i}([Z,X_{i}]),&[Z,1\otimes
  1]>=<[Z^{*},\delta_{i}([Z,X_{i}])],1\otimes  1>\\
  &=<\delta_{i}([Z^{*},[Z,X_{i}]]),1\otimes
  1>-<[\delta_{i}(Z^{*}),[Z,X_{i}]],1\otimes  1>\\
  &=<[Z,X_{i}],[Z,\Delta(X_{i})]>+<\delta_{i}(Z^{*}),[1\otimes
  1,[Z^{*},X_{i}]]>\end{align*}
We can now conclude using that $1\otimes\tau(\delta_{i}(Z^{*}))=\tau\otimes1(\delta_{i}(Z))^{*}$: 
\begin{align*}n &||[Z,1\otimes  1]||^{2}_{HS}=||[Z,1\otimes
1]||^{2}\\&+\sum_{i=1}^{n}<[Z,X_{i}],[Z,\Delta(X_{i})]>+2\sum_{i=1}^{n}\Re<\delta_{i}(Z),[1\otimes1,[Z,X_{i}]]>.\end{align*}\end{preuves}

%Recall from Theorem 1 in \citep{S05} the characterization of finite Fisher information in terms of dual systems of unbounded operators, i.e. (with the notation of the lemma), 
%$\Phi^{*}(X_{1},...,X_{n})<\infty$ if and only if there exists closable unbounded operators $Y_{i}:L^{2}(W^{*}(X_{1},...,X_{n}))\rightarrow L^{2}(W^{*}(X_{1},...,X_{n}))$ whose domains includes $\C\langle X_{1},...,X_{n} \rangle$ so that $Y_{i}1=0$ and 1 belongs to the domain of $Y_{i}^{*}$ and so that $[Y_{j},X_{i}]=\delta_{j}(X_{i})$.
%In that case, $\delta_{j}$ clearly coincide with $[Y_{j},.]$ on $\C\langle X_{1},...,X_{n} \rangle$ (with those derivations seen as valued in Hilbert-Schmidt operators), and moreover as a consequence of the proof in\citep{S05}, $\delta_{j}^{*}(1\otimes 1)=JY_{j}^{*}1$ and $\C\langle X_{1},...,X_{n} \rangle$ is in the domain of $Y_{j}^{*}$ ($j=1..n$).

%The following lemma use these results to reformulate the previous lemma.
%\begin{lemma}\label{FPC}
%With the same assumptions as in the previous lemma :$$(n-1)||Z-\tau(Z)||_{2}^{2}=\sum_{i=1}^{n} <[Z,[Z,X_{i}]],\Delta(X_i)>+<[Z,X_{i}],\overline{\delta_{i}}(Z).1+J\overline{\delta_{i}(Z)}.1>.$$
%\end{lemma}

%\begin{preuves}
%Since $\C\langle X_{1},...,X_{n} \rangle$ is a core for $\overline{\delta}$, it is sufficient to prove this for $Z\in \C\langle X_{1},...,X_{n} \rangle$ and self-adjoint. Moreover, since $1\otimes\tau\circ\delta_{i}$ is easily seen to be closable like $\delta_{i}$, this inequality is even valid for $Z\in \cap_{i=1}^{n}D(\overline{1\otimes\tau\circ\delta_{i}})$
%\end{preuves}

For our purpose, the following lemma is only an intermediary step to the next lemma, but, as the remark after it shows, it can have an independent interest.
\begin{lemma}\label{FPC2}
Let $(M,\tau)$ a tracial $W^{*}$-probability space. Let $(X_{1},...X_{n})$
a n-tuple of $n\geq 2$ self-adjoints such that the microstates free Fisher information
$\Phi^{*}(X_{1},...,X_{n})<\infty$. Let $Z \in W^{*}(X_{1},...,X_{n})\cap
D(\bar{\delta})$, then
 the following inequality holds :
$$||(1\otimes\tau)(\bar{\delta_{i}}(Z))-Z\Delta(X_{i})||_{2}^{2}\leq ||Z\Delta(X_{i})||_{2}^{2}+ | \langle\bar{\delta_{i}}(Z^{*}Z),1\otimes\Delta(X_{i})\rangle|$$
Thus, if we assume moreover we have second order conjugate variables $\I_{2,j}=\I_{2}(X_{j}:C\langle X_{1},...,\hat{X_{j}},...,X_{n}\rangle)$ (in $L^{1}(M,\tau)$, as defined in \citep{Vo5} Definition 3.1). Then
 the following inequality holds :
$$||(1\otimes\tau)(\bar{\delta_{i}}(Z))||_{2}\leq 2||Z\Delta(X_{i})||_{2}+  | \langle Z^{*}Z,\I_{2,i}\rangle|^{1/2}.$$
As a consequence, $(1\otimes\tau)\circ\bar{\delta_{i}}$ extends as a bounded map $M\rightarrow L^{2}(M,\tau)$.% and :
%\begin{align*}2(n-1)&||Z-\tau(Z)||_{2}^{2}\leq\sum_{i=1}^{n} ||[Z,X_{i}]||_{2}\left(9||Z\Delta(X_{i})||_{2}+||\Delta(X_{i})Z||_{2}+6| \langle Z^{2},\I_{2,i}\rangle|^{1/2}\right).\end{align*}
\end{lemma}

\begin{remark}
If we assume moreover we have bounded first and second order conjugate variables (i.e. $\Delta(X_{i}),\I_{2,i}\in M$ or more generally bounded conjugate variable and dual system $\tilde{Y}_{i}$ in the sense of Voiculescu), then the previous lemma shows that $(1\otimes\tau)\circ\bar{\delta_{i}}$ extends as a bounded map on $L^{2}(M)$and moreover the inequality above implies that  $X_{1},...,X_{n}$ is a non-$\Gamma$ set (for $D(\bar{\delta})$) in the sense of \citep{P2}. As a consequence of Corollary 3.3 in \citep{P2} this shows that $W^{*}(X_{1},...,X_{n})$ doesn't have property (T). A study under less restrictive assumptions will need new investigations, but we can already note (using well-known results of \citep{Vo5}) that this implies that for any $X_{1},...,X_{n}$, if $S_{1},...,S_{n}$ is a free semicircular system free with $X_{1},...X_{n}$, then $W^{*}(X_{1}+\epsilon S_{1},...,X_{n}+\epsilon S_{n})$ doesn't have property (T) (for any $\epsilon>0$). This result was proved in \citep{SJ} assuming moreover $W^{*}(X_{1},...,X_{n})$ embeddable in $R^{\omega}$.
\end{remark}

\begin{preuves}
The only non-trivial statement is the first one (using the previous lemma for proving consequences). Moreover we can assume $Z\in \C\langle X_{1},...,X_{n}\rangle$ as usual (take the limit in the fifth line below using proposition \ref{Wdev} and then compute with $Z\in D(\bar{\delta})\cap M$). Let us compute (using the formula for $\delta_{i}^{*}$, Corollary 4.3 in \citep{Vo5}, and coassociativity in the third line):

\begin{align*}||(1\otimes\tau)&(\bar{\delta_{i}}(Z))||_{2}^{2} =\langle (1\otimes\tau)(\bar{\delta_{i}}(Z))\otimes 1,\bar{\delta_{i}}(Z)\rangle\\ 
&=\langle (1\otimes\tau)(\bar{\delta_{i}}(Z))\Delta(X_{i}),Z\rangle -\langle(1\otimes\tau)\bar{\delta_{i}}(1\otimes\tau)(\bar{\delta_{i}}(Z)),Z\rangle\\
&=\langle (1\otimes\tau)(\bar{\delta_{i}}(Z))\Delta(X_{i}),Z\rangle -\langle(1\otimes\tau\otimes \tau)1\otimes\bar{\delta_{i}}\circ\bar{\delta_{i}}(Z)),Z\rangle\\
&=\langle (1\otimes\tau)(\bar{\delta_{i}}(Z))\Delta(X_{i})-(1\otimes\tau)(\bar{\delta_{i}}(Z)\Delta(X_{i})),Z\rangle \\ 
&=\langle (\bar{\delta_{i}}(Z)),Z[\Delta(X_{i}),1\otimes 1]\rangle \\ 
&=\langle (\bar{\delta_{i}}(Z)),Z\Delta(X_{i})\otimes 1\rangle -\langle (\bar{\delta_{i}}(Z^{*}Z)-\bar{\delta_{i}}(Z^{*})Z,1\otimes \Delta(X_{i})\rangle\\ 
&=\langle (\bar{\delta_{i}}(Z)),Z\Delta(X_{i})\otimes 1\rangle+\langle (\bar{\delta_{i}}(Z^{*})),1\otimes \Delta(X_{i})Z^{*}\rangle \\ &\ \ \ \ \ -\langle (\bar{\delta_{i}}(Z^{*}Z)),1\otimes \Delta(X_{i})\rangle. 
%\\ &\leq ||(1\otimes\tau)(\bar{\delta_{i}}(Z))||_{2}(||Z\Delta(X_{i})||_{2}+||\Delta(X_{i})Z^{*}||_{2})+\left|\langle (\bar{\delta_{i}}(Z^{*}Z)),1\otimes \Delta(X_{i})\rangle\right|. 
\end{align*}

Now note we can use $(1\otimes\tau)(\bar{\delta_{i}}(Z))^{*}=(\tau\otimes 1)(\bar{\delta_{i}}(Z^{*}))$ (using $\bar{\delta_{i}}$ is a real derivation), to conclude :

\begin{align*}||(1\otimes\tau)(\bar{\delta_{i}}(Z))-Z\Delta(X_{i})||_{2}^{2}
&=||Z\Delta(X_{i})||_{2}^{2} -\langle (\bar{\delta_{i}}(Z^{*}Z)),1\otimes \Delta(X_{i})\rangle. 
%\\ &\leq ||(1\otimes\tau)(\bar{\delta_{i}}(Z))||_{2}(||Z\Delta(X_{i})||_{2}+||\Delta(X_{i})Z^{*}||_{2})+\left|\langle (\bar{\delta_{i}}(Z^{*}Z)),1\otimes \Delta(X_{i})\rangle\right|. 
\end{align*}

%Now, a usual inequality associated to a degree two polynomial concludes.
%Likewise, we get  :

%\begin{align*}||(\tau\otimes1)(\bar{\delta_{i}}(Z))||_{2}^{2}&=\langle 1\otimes(\tau\otimes1)(\bar{\delta_{i}}(Z)),\bar{\delta_{i}}(Z)\rangle\\ 
%&=\langle \Delta(X_{i})(\tau\otimes1)(\bar{\delta_{i}}(Z)),Z\rangle -\langle(\tau\otimes1)\bar{\delta_{i}}(\tau\otimes1)(\bar{\delta_{i}}(Z)),Z\rangle\\
%&=\langle \Delta(X_{i})(\tau\otimes1)(\bar{\delta_{i}}(Z)),Z\rangle -\langle(\tau\otimes \tau\otimes1)1\otimes\bar{\delta_{i}}\circ\bar{\delta_{i}}(Z)),Z\rangle\\
%&=\langle \Delta(X_{i})(\tau\otimes1)(\bar{\delta_{i}}(Z))-(\tau\otimes1)(\Delta(X_{i})\bar{\delta_{i}}(Z)),Z\rangle \\ 
%&=\langle (\bar{\delta_{i}}(Z)),[1\otimes 1,\Delta(X_{i})]Z\rangle \\ 
%&=\langle (\bar{\delta_{i}}(Z)),1\otimes\Delta(X_{i})Z\rangle -\langle (\bar{\delta_{i}}(Z}^{2))-Z(\bar{\delta_{i}}(Z}), \Delta(X_{i})\otimes1\rangle\\ 
%&=\langle (\bar{\delta_{i}}(Z)),Z\Delta(X_{i})\otimes 1+1\otimes \Delta(X_{i})Z\rangle -\langle (\bar{\delta_{i}}(Z}^{2)),\Delta(X_{i})\otimes 1\rangle. 
%\end{align*}
\end{preuves}

\begin{lemma}\label{FPC3}
Let $(M,\tau)$ a tracial $W^{*}$-probability space. Let $(X_{1},...X_{n})$
a n-tuple of  $n\geq 2$ self-adjoints such that the microstates free Fisher information
$\Phi^{*}(X_{1},...,X_{n})<\infty$. Let $Z \in W^{*}(X_{1},...,X_{n})\cap
D(\bar{\delta})$ a selfadjoint or a unitary, then
we have the following inequality :
$$||(1\otimes\tau)(\bar{\delta_{i}}(Z))-Z\Delta(X_{i})||_{2}\leq ||\Delta(X_{i})||_{2}||Z||$$
As a consequence, $(1\otimes\tau)\circ\bar{\delta_{i}}$ extends as a bounded map $M\rightarrow L^{2}(M,\tau)$ and :
\begin{align*}(n-1)&||Z-\tau(Z)||_{2}^{2}\leq\sum_{i=1}^{n} -\frac{1}{2}<[Z,X_{i}],[Z,\Delta(X_{i})]>+2||[Z,X_{i}]||_{2}||\Delta(X_{i})||_{2}||Z||.\end{align*}
\end{lemma}

\begin{preuves}
Take $Z$ of norm less than 1 ($||Z||<1$).
If it is selfadjoint, we can write $Z$ as an half sum of two unitaries in $W^{*}(X_{1},...,X_{n})\cap
D(\bar{\delta})$ using stability by $C^{1}$ functional calculus (e.g. lemma 7.2 in \citep{CiS}), we have only to prove the inequality for any unitary $U$. 
This follows at once from the previous lemma. The second statement is a direct consequence.
\end{preuves}

\subsection{The main result}

Now lemma \ref{FPC3} contains immediately the non-$\Gamma$ result we wanted :

\begin{theorem}
Let $(M,\tau)$ a tracial $W^{*}$-probability space. Let $(X_{1},...X_{n})$
a n-tuple (of self-adjoints, $n\geq 2$) such that the microstates free Fisher information
$\Phi^{*}(X_{1},...,X_{n})<\infty$, then $W^{*}(X_{1},...,X_{n})$
doesn't have property $\Gamma$, i.e. all central sequences $Z_{m}$ (i.e bounded in
$W^{*}(X_{1},...,X_{n})$  and such that $\forall Y \in
W^{*}(X_{1},...,X_{n}) ||[Z_{m}, Y]||_{2}\rightarrow 0$)
 are trivial : $||Z_{m} -\tau(Z_{m})||_{2}\rightarrow 0$. As a consequence,  $W^{*}(X_{1},...,X_{n})$ is not amenable.
\end{theorem}

\section{Factoriality under finite non-microstates entropy}

We will now prove factoriality under the weaker assumption $\chi^{*}(X_{1},...,X_{n})>-\infty$ as another consequence of lemma \ref{FPC3}.

In this part, we thus let $X_{1},...,X_{n}$ $n\geq 2$ selfadjoints, $S_{1},...,S_{n}$, a free semi-circular system, free with $X_{1},...,X_{n}$. Let $Y_{j}^{t}=X_{j}+\sqrt{t}S_{j}$ and $E_{t}$ the (trace preserving) conditional expectation onto $W^{*}(Y_{1}^{t},...,Y_{n}^{t})$(seen as a sub-von Neumann algebra of $W^{*}(\{X_{i},S_{i}\})$). Then recall that the non-microstates entropy is defined by the following integral :

\begin{align*}\chi^{*}&(X_{1},...,X_{n})=\frac{1}{2}\int_{0}^{\infty}\left(\frac{n}{1+t}-\Phi^{*}(X_{1}+\sqrt{t}S_{1},...,X_{n}+\sqrt{t}S_{n})\right)dt +\frac{n}{2}\log 2\pi e.\end{align*}

It is readily seen that if $\liminf_{t\rightarrow 0} t\Phi^{*}(X_{1}+\sqrt{t}S_{1},...,X_{n}+\sqrt{t}S_{n})\neq 0$ we have necessarily $\chi^{*}(X_{1},...,X_{n})=-\infty$, we will thus use the assumption in that way. More generally, a variant of free entropy dimension was defined in \citep{CS} by $\delta^{\star}(X_{1},...,X_{n})=n-\liminf_{t\rightarrow 0} t\Phi^{*}(X_{1}+\sqrt{t}S_{1},...,X_{n}+\sqrt{t}S_{n})$, we will thus express our result in function of this entropy dimension. %Moreover, we will use the fact that we know from \citep{Vo5} (Corollary 3.9) the second order conjugate variable for $Y_{1}^{t},...,Y_{n}^{t}$ and benefit from this via lemma \ref{FPC2}:
%$$\I_{2,j}^{t}=\I_{2}(Y_{j}^{t}:C\langle Y_{1}^{t},...,\hat{Y_{j}^{t}},...,Y_{n}^{t}\rangle)=t^{-1}E_{t}(S_{j}^{2}-1).$$

We can now prove our claimed result :

\begin{theorem}
Let $(M,\tau)$ a tracial $W^{*}$-probability space. Let $(X_{1},...X_{n})$
a n-tuple (of self-adjoints, $n\geq 2$) then the following inequality holds for any central selfadjoint $Z$ :
$$||Z-\tau(Z)||_{2}^{2}\leq 2\frac{n- \delta^{\star}(X_{1},...,X_{n})}{n-1}||Z||^{2}.$$
As a consequence, 
if $n- \delta^{*}(X_{1},...,X_{n}):=\liminf_{t\rightarrow 0} t\Phi^{*}(X_{1}+\sqrt{t}S_{1},...,X_{n}+\sqrt{t}S_{n})=0$, then $W^{*}(X_{1},...,X_{n})$ is a factor. As another example, if $\delta^{*}(X_{1},...,X_{n})> \frac{n+1}{2}$, $W^{*}(X_{1},...,X_{n})$ has no central projection of trace one half, especially, doesn't have diffuse center.
\end{theorem}

\begin{preuves}
Let $Z$ in the center of $W^{*}(X_{1},...,X_{n})$, with $||Z||\leq 1$ and apply lemma \ref{FPC3} to $E_{t}(Z)\in W^{*}(Y_{1}^{t},...,Y_{n}^{t})$ to get~:
\begin{align*}(&n-1)||E_{t}(Z)-\tau(E_{t}(Z))||_{2}^{2}\\ &\leq\sum_{i=1}^{n} -\frac{1}{2}<[E_{t}(Z),Y_{i}^{t}],[E_{t}(Z),\Delta(Y_{i}^{t})]>+2||[E_{t}(Z),Y_{i}^{t}]||_{2}||\Delta(Y_{i}^{t})||_{2}||Z||\\ &=\sum_{i=1}^{n} -\frac{1}{2}<E_{t}([Z,Y_{i}^{t}]),[E_{t}(Z),\frac{1}{\sqrt{t}}E_{t}(S_{i})]>+2||E_{t}([Z,Y_{i}^{t}])||_{2}||\frac{1}{\sqrt{t}}E_{t}(S_{i})||_{2}||Z||\\ &=\sum_{i=1}^{n} -\frac{1}{2}<E_{t}([Z,S_{i}]),[E_{t}(Z),E_{t}(S_{i})]>+2||E_{t}([Z,S_{i}])||_{2}||E_{t}(S_{i})||_{2}||Z||\end{align*}\begin{align*}(&n-1)||E_{t}(Z)-\tau(E_{t}(Z))||_{2}^{2}\\ &\leq\sum_{i
=1}^{n} -\frac{1}{2}<E_{t}([Z-E_{t}(Z),S_{i}]),[E_{t}(Z),E_{t}(S_{i})]> -\frac{1}{2}<[E_{t}(Z),E_{t}(S_{i})],[E_{t}(Z),E_{t}(S_{i})]>\\ &\ \ \ \ \ \ +2(||E_{t}([Z-E_{t}(Z),S_{i}])||_{2}+||[E_{t}(Z),E_{t}(S_{i})]||_{2})||E_{t}(S_{i})||_{2}||Z||\\ &\leq\sum_{i=1}^{n} 12||Z-E_{t}(Z)||_{2}||E_{t}(S_{i})||_{2}||Z||-\frac{1}{2}||[E_{t}(Z),E_{t}(S_{i})]||_{2}^{2}\\&\ \ \ \ \ +||[E_{t}(Z),E_{t}(S_{i})]||_{2}^{2}+(||E_{t}(S_{i})||_{2}||Z||)^{2}\\ &\leq 12n||Z-E_{t}(Z)||_{2}||Z||+\sum_{i=1}^{n}\frac{1}{2}||[E_{t}(Z),E_{t}(S_{i})]||_{2}^{2}+(||E_{t}(S_{i})||_{2}||Z||)^{2}\\ &\leq 12n||Z-E_{t}(Z)||_{2}||Z||+2||Z||^{2}\sum_{i=1}^{n}||E_{t}(S_{i})||_{2}^{2}.\end{align*}

We used at the second line the result of \citep{Vo5} about the conjugate variable in the algebra generated by $Y_{i}^{t}$ :$\Delta(Y_{i}^{t})=\frac{1}{\sqrt{t}}E_{t}(S_{i})$. We also used conditional expectation property and then at line 3 that $Z$ commutes with $X_{i}$. In the fourth line we used $Z=Z-E_{t}(Z)+E_{t}(Z)$ and then we only compute using $||S_{i}||=2$, $||S_{i}||_{2}=1$ and arithmetico geometric inequality.

At the end, we thus get using the definition of free Fisher information and the result of \citep{Vo5} above:
$$(n-1)||E_{t}(Z)-\tau(E_{t}(Z))||_{2}^{2}\leq 12n||Z-E_{t}(Z)||_{2}||Z||+2||Z||^{2}t\Phi^{*}(X_{1}+\sqrt{t}S_{1},...,X_{n}+\sqrt{t}S_{n}).$$

It is thus sufficient to notice that $||E_{t}(Z)-Z||_{2}$ goes to 0 with $t$ to get the inequality stated by taking a $\liminf$. $||E_{t}(Z)-Z||_{2}\rightarrow 0$ follows from Kaplansky density theorem, and from the remark that for $P$ a non commutative polynomial $||E_{t}(P(X_{1},...,X_{n}))-P(X_{1},...,X_{n})||_{2}\leq ||P(Y_{1}^{t},...,Y_{n}^{t})-P(X_{1},...,X_{n})||_{2}$.

The consequences are trivial : for instance for the second, apply the inequality to $Z=1-2P$ the corresponding central selfadjoint unitary of trace $0$ if $P$ a central projection of trace $1/2$.
\end{preuves}

\textbf{Acknowledgments}
The author would like to thank Professor Dan Voiculescu for allowing him to use the argument of his free Poincaré inequality in another context, Professor Dimitri Shlyakhtenko for showing
him references \citep{S3}, \citep{CiS} and \citep{P1}, and Professor Jesse Peterson for pointing out to his attention his result Corollary 3.3 in \citep{P2}. The author would also
like to thank D. Shlyakhtenko, J. Peterson and P. Biane for useful
discussions and comments.

\bibliographystyle{plain}
\bibliography{microstate-free-non-Gamma-5}

\end{document}